\documentstyle[12pt]{article}
\title{\bf Cauchy-Riemann Geometry and Contact Topology
     in Three Dimensions}
\author{
  Jih-Hsin Cheng \\
Institute of Mathematics \\
Academia Sinica \\
Taipei, Taiwan, R.O.C.}
\date{ }
\newenvironment{indention}[1]{\par
\addtolength{\leftskip}{#1}
\begingroup}{\endgroup\par}

\begin{document}
\maketitle
{\large \setlength{\baselineskip}{24pt}
\begin{abstract}
We introduce a global Cauchy-Riemann($CR$)-invariant and discuss its
behavior on
the moduli space of $CR$-structures. We argue that this study
is related to the Smale conjecture in 3-topology and the problem of counting
complex structures. Furthermore, we propose a contact-analogue of
Ray-Singer's analytic torsion. This
``contact torsion'' is expected to be able to distinguish among ``contact
lens'' spaces. We also propose the study of a certain kind of monopole
equation associated with a contact structure.
\end{abstract}



\section{Introduction}

We study low-dimensional problems in topology and geometry via
a study of contact and Cauchy-Riemann ($CR$) structures. Let us start
with a closed (compact without boundary) oriented three-manifold
$M$. A contact structure (or bundle) $\xi$ on $M$ is a completely
non-integrable rank 2 subbundle of $TM$. It is well known that
there are no local invariants for contact structures according to a 
classical theorem of Darboux. Also, two nearby contact structures
on a closed manifold are isotopy-equivalent by Gray's theorem
(Gray, 1959; Hamilton, 1982). Therefore, a contact structure has
no continuous moduli. In 
this sense, it is a kind of geometric structure even softer than a
complex structure. The isotopy classes are distinguished by so-called
tight or overtwisted contact structures (Eliashberg, 1992). The existence of
contact structures on a closed oriented three-manifold is known from
the work of Martinet (1971) and Lutz (1977). (See also Altschuler (1995)
for an analytic proof using the so-called linear contact flow.)

Given a contact structure, we can consider a $CR$-structure,
i.e., a ``complex structure'' defined on a contact bundle.
Different from the usual complex structure, a $CR$-structure does
have local invariants. Thus, analysis is needed. In Section II, we give a
brief introduction to $CR$-geometry and an application in K$\ddot{a}$hler
geometry. In Section III, we introduce a global $CR$-invariant $\mu_{\xi}$
and discuss
its behavior on the moduli space of $CR$-structures. Also, we argue
that the contractibility of our $CR$ moduli space for $S^3$ 
confirms the so-called Smale conjecture.

In Section IV, we discuss spherical
$CR$-structures: the critical points of $\mu_{\xi}$. To distinguish
among ``$CR$ lens'' spaces, we propose a possible $CR$-invariant defined
for spherical $CR$-structures, which is a contact-analogue of
Ray-Singer's analytic torsion. In Section V, we give a heuristic argument
for how our understanding of $\mu_{\xi}$ can be applied to the problem
of counting the number of complex structures on a closed four-manifold.
In Section VI, we propose the study of a certain kind of monopole equation 
for contact three-manifolds. 

   
\bigskip

\section{Basics in $CR$-Geometry}

A $CR$-structure $J$ compatible with the contact structure $\xi$ is 
a complex structure on $\xi$, i.e., a bundle endomorphism $J:{\xi}
{\rightarrow}{\xi}$, such that $J^{2}=-Identity$. Natural examples
come from boundaries of strictly pseudoconvex domains $D$ in $C^2$.
Let $J_{C^2}$ denote the multiplication by $i$ in $C^2$. Let our
three-manifold $M={\partial}D$, the boundary of $D$. The contact
structure $\xi$ is considered to be the intersection of $TM$ and
$J_{C^2}TM$, the tangent subspaces invariant under $J_{C^2}$.
In addition, our $CR$-structure is taken to be a restriction of $J_{C^2}$
on $\xi$. This $CR$-structure is usually called the $CR$-structure
induced from $C^2$.

In his famous theorem, Fefferman (1974) asserts that two strictly
pseudoconvex domains with smooth boundaries in $C^{n+1}$ are 
biholomorphic to each other if and only if their boundaries are
$CR$-equivalent. Therefore the $CR$-structure on the boundary reflects 
the complex structure of the inside domain. It is well known that
we have the Riemann mapping theorem in $C^1$. However, this theorem
is no longer true for higher dimensions. Indeed, we do have local
invariants for our $CR$ manifold $(M, {\xi}, J)$ (e.g., Cartan (1932)
and Chern and Moser (1974)).
 
First, choose eigenvectors $Z_1$, $Z_{\bar 1}$ of $J$ with eigenvalues
$i$, $-i$, respectively. Let $\{ {\theta^1},{\theta^{\bar 1}} \}$ be a set
of complex one-forms dual to $\{ Z_{1}, Z_{\bar 1} \}$. Then,
choose a local one-form $\theta$ annihilating $\xi$ (called contact
form) so that

$$
d{\theta}=ih_{1{\bar 1}}{\theta^1}{\wedge}{\theta^{\bar 1}}+
          {\theta}{\wedge}{\phi}
$$

\noindent for some real one-form $\phi$ and positive $h_{1{\bar 1}}$.
(We will use $h_{1{\bar 1}}$ and $h^{1{\bar 1}}=(h_{1{\bar 1}})^{-1}$
to raise or lower indices.) Now, for a different choice of coframe
$({\tilde{\theta}}, {\tilde{\theta}}^1, {\tilde{\theta}}^{\bar 1};
{\tilde{\phi}})$ satisfying the above equation, we have the following
transformation relation:

\begin{equation}
\left\{ \begin{array}{lll}
        {\tilde{\theta}}=u{\theta}\\
        {\tilde{\theta}}^1={u_1}^{1}{\theta}^1+v^{1}{\theta}\\
        {\tilde{\phi}}=-\frac{du}{u}+{\phi}+
           2Re(iu^{-1}v^{1}u_{{\bar 1}1}{\theta}^{\bar 1})+s{\theta}
        \end{array}
\right.
\end{equation}

\noindent for positive $u$ and some real function $s$. Differentiating
${\theta}^1$, $\phi$ gives the first structural equations:

\begin{equation}
\left\{ \begin{array}{ll}
  d{\theta}^{1}={\theta}^{1}{\wedge}{\phi_1}^{1}+{\theta}{\wedge}{\phi^{1}}\\
  d{\phi}=2Re(i{\theta_{\bar 1}}{\wedge}{\phi^{\bar 1}})+
         {\theta}{\wedge}{\psi} 
        \end{array} 
\right.
\end{equation}

\noindent for the connection forms ${\phi_1}^{1},{\phi^{1}},{\psi}$.
Differentiating the connection forms again and requiring certain trace
conditions (e.g., Chern and Moser (1974)), we obtain the second set
of structural equations:

\begin{equation}
\left\{ \begin{array}{lll}
  d{\phi_1}^{1}-i{\theta_{1}}{\wedge}{\phi^{1}}+2i{\phi_1}{\wedge}{\theta^1}
       +\frac{1}{2}{\psi}{\wedge}{\theta}=0   \\
  d{\phi^{1}}- {\phi}{\wedge}{\phi^{1}}-{\phi^{1}}{\wedge}{\phi_1}^{1}
  +\frac{1}{2}{\psi}{\wedge}{\theta^1}={Q^1}_{\bar 1}{\theta^{\bar 1}}
   {\wedge}{\theta}   \\
  d{\psi}-{\phi}{\wedge}{\psi}-2i{\phi^{1}}{\wedge}{\phi_{1}}=
      (R_{1}{\theta^1}+R_{\bar 1}{\theta^{\bar 1}}){\wedge}{\theta}, 
 \end{array}
\right.
\end{equation}

\noindent in which ${Q^1}_{\bar 1}$ or $Q_{11}$ is called the Cartan
(curvature) tensor, and $R_{1}, R_{\bar 1}$ are determined by means of
suitable covariant derivatives of ${Q^1}_{\bar 1}$ (Cheng, 1987).
The normalization condition: ${\phi}-{\phi_1}^{1}-{\phi_{\bar 1}}^{\bar 1}
=0$ and the above structural equations Eqs.(2) and (3) uniquely
determine the connection forms ${\phi_1}^{1},{\phi^{1}},{\psi}$. 
Under the change of coframe Eq.(1), the Cartan tensor is transformed as
follows:

$$ Q_{11}={\tilde Q}_{11}u({u_1}^{1})^{2}. $$

The fundamental theorem of 3-dimensional $CR$-geometry due to
Cartan (1932a,1932b) asserts that $Q_{11}=0$ if and only if $(M,{\xi},J)$
is locally $CR$-equivalent to $(S^{3},{\hat{\xi}},{\hat{J}})$, where
$({\hat{\xi}},{\hat{J}})$ denotes the standard $CR$-structure on the
unit 3-sphere $S^3$, induced from $C^2$.

\bigskip
\noindent $\bf Definition$. We call a $CR$ manifold $(M,{\xi},J)$ or
just $J$ spherical if it is locally $CR$-equivalent to 
$(S^{3},{\hat{\xi}},{\hat{J}})$. Quantitatively, a $CR$-structure is
spherical if $Q_{11}=0$ according to Cartan's theorem. 

\bigskip
\noindent $\bf{An\:Application\:in\:K\ddot{a}hler\:Geometry}$

\medskip
Let $N$ be an $n$-dimensional K$\ddot{a}$hler manifold. Suppose we have
a holomorphic line bundle $L$ with the first Chern class being the
K$\ddot{a}$hler class so that a suitable circle bundle $M{\subset}
L$ with the induced $CR$-structure is closely related to the K$\ddot{a}$hler
geometry of $N$ (Webster, 1977). It turns out that we can identify (up to
a constant) the Cartan
tensor $Q_{11}$ of $M$ with $R_{,11}$, the covariant derivative of the scalar
curvature $R$ of $N$ in the (1,0)-direction twice. When $n{\geq}2$,
we can identify the Chern tensor (Chern and Moser, 1974, 1983) in
higher dimensional
$CR$-geometry with the Bochner tensor of $N$. In 1977, Sid Webster
applied $CR$-geometry to obtain the following result:

\begin{indention}{2ex}
\noindent
\normalsize
\\ 
Let $N$ be a simply-connected closed
K$\ddot{a}$hler manifold of dimension $n$. Suppose that $N$ admits a Hodge 
metric for which the Bochner tensor vanishes
if $n{\geq}2$ or for which $R_{,11}$ vanishes if $n=1$. Then, $N$
is holomorphically
isometric to complex projective space $CP^n$ with a standard
Fubini-Study metric. (Webster, 1977)
\\
\end{indention}

Next, relative to a special coframe $({\theta},{\theta^1},{\theta^{\bar 1}};
{\phi}=0)$ satisfying $d{\theta}=ih_{1{\bar 1}}{\theta^1}{\wedge}
{\theta^{\bar 1}}$, we can define the so-called pseudohermitian
connection ${\omega_1}^1$, torsion $A_{11}$, and curvature $\cal W$,
called the Tanaka-Webster curvature (Tanaka, 1975; Webster, 1977).
These data are uniquely determined
by the following equations:

$$ \left
 \{ \begin{array}{lll}
   d{\theta^1}={\theta^1}{\wedge}{\omega_1}^1+{A^1}_{\bar 1}
               {\theta}{\wedge}{\theta^{\bar 1}}\\
   d{\omega_1}^1={\cal W}{\theta^1}{\wedge}{\theta^{\bar 1}}
               \:(mod\:{\theta})\\
   {\omega_1}^1+{\omega_{\bar 1}}^{\bar 1}=h^{1{\bar 1}}dh_{1{\bar 1}}.
   \end{array} \right.
$$ 

The torsion $A_{11}$ and the Tanaka-Webster curvature $\cal W$ are not
``tensorial'' under the change of contact form ${\tilde{\theta}}=
u{\theta}$ (Lee, 1986), but are ``tensorial'' under the change
${\tilde{\theta}}^
{1}={u_1}^{1}{\theta^1}$. The Cartan tensor can be expressed in terms
of these data (Cheng and Lee, 1990):

$$
Q_{11}=\frac{1}{6}{\cal W}_{,11}+\frac{i}{2}{\cal W}A_{11}-A_{11,0}
       -\frac{2i}{3}A_{11,{\bar 1}1}.
$$   

\noindent Here, covariant derivatives are taken with respect to the
pseudohermitian connection ${\omega_1}^1$, and ``0'' means the $T$-
direction. (The tangent vector field $T$ is uniquely determined by
${\theta}(T)=1$ and $L_{T}{\theta}=0$.) Before going on, 
another result should be noted:

\begin{indention}{2ex}
\noindent
\normalsize
\\
The boundary of a circular
domain in $C^{n+1}$ is $CR$-equivalent to the unit sphere $S^{2n+1}
{\subset}C^{n+1}$ with the standard induced $CR$-structure if and only
if the Tanaka-Webster curvature ${\cal W}{\equiv}constant$ (with respect
to a suitable choice of contact form)(Unpublished paper by J. Bland and
P. M. Wang).
\\
\end{indention}

The proof of the above result in the original draft contains a gap which
can be remedied by the following result:

\begin{indention}{2ex}
\noindent
\normalsize
\\
Let $N$ be a closed complex
manifold with two K$\ddot{a}$hler metrics $g$, $\tilde{g}$. Suppose
the Bochner tensor of $g$ vanishes and the scalar curvature of $\tilde{g}$
is a constant. Then, the fact that the K$\ddot{a}$hler class of $g$ is
cohomologous
to the K$\ddot{a}$hler class of $\tilde{g}$ implies that $(N,g)$ and
$(N,{\tilde g})$ are isometric to each other (Chen and Lue, 1981).
\\
\end{indention}
\bigskip

\section{The $\mu_{\xi}$-Invariant and the Moduli Space}

First, we will construct an energy functional on
the space of $CR$-structures so that the critical points consist
of spherical $CR$-structures. Let $\Pi$ denote the $su(2,1)$-valued
Cartan connection form defined by

$$
{\Pi}=\left ( \begin{array}{ccc}
      -\frac{1}{3}({\phi_1}^{1}+{\phi}) & {\theta^1} & 2{\theta}\\
      -i{\phi_1} & \frac{1}{3}(2{\phi_1}^{1}-{\phi}) & 2i{\theta_1}\\
      -\frac{1}{4}{\psi} & \frac{1}{2}{\phi^1} & \frac{1}{3}({\phi}+
      {\phi_{\bar 1}}^{\bar 1}) \end{array} \right ).
$$

The curvature form $\Omega$ is defined as usual by ${\Omega}=d{\Pi}-
{\Pi}{\wedge}{\Pi}$. The transgression $TC_{2}({\Pi})$ of the second
Chern form is given by

\begin{eqnarray}
TC_{2}({\Pi})&=&\frac{1}{8{\pi^2}}[tr({\Pi}{\wedge}{\Omega}+\frac{1}{3}
       tr({\Pi}{\wedge}{\Pi}{\wedge}{\Pi})]\nonumber\\
  &=&\frac{1}{24{\pi^2}}tr({\Pi}{\wedge}{\Pi}{\wedge}{\Pi})\nonumber\\  
  & &(since\:tr({\Pi}{\wedge}{\Omega})=0).\nonumber
\end{eqnarray}

We can verify that
the 3-form $TC_{2}({\Pi})$ is invariant under the
change of contact form and invariant up to an exact form under the 
coframe change Eq.(1). In the late 1980's, Burns and Epstein (1988)(also
Cheng and Lee (1990)) discovered that the integral of
$TC_{2}({\Pi})$,
denoted as $\mu_{\xi}$,
is a global $CR$-invariant (assuming trivial holomorphic tangent
bundle as in Burns and Epstein (1988); extended to arbitrary $M$
by a relative version of the invariant in Cheng and Lee (1990)): 

\medskip
\begin{eqnarray}
{\mu_{\xi}}(J)&=&\frac{1}{24{\pi^2}}{\int_M}tr
                 ({\Pi}{\wedge}{\Pi}{\wedge}{\Pi})\nonumber\\
 &=&\frac{1}{8{\pi^2}}{\int_M}[2Re(i{\theta^1}{\wedge}{\phi^{\bar 1}}{\wedge}
   {\phi_1}^{1})+\frac{1}{2}{\theta}{\wedge}{\psi}{\wedge}{\phi}
    -2i{\theta}{\wedge}{\phi^1}{\wedge}{\phi^{\bar 1}}
       -\frac{1}{2}d({\theta}{\wedge}{\psi})]\nonumber\\
 &=&\frac{1}{8{\pi^2}}{\int_M}[(\frac{1}{6}{\cal W}^{2}+2|A_{11}|^{2})
   {\theta}{\wedge}d{\theta}+\frac{2}{3}{\omega_1}^{1}{\wedge}d{\omega_1}^{1}]
   \nonumber\\
 & &(in\:terms\:of\:pseudohermitian\:geometry).\nonumber
\end{eqnarray}

\medskip
It is remarkable that the above integral is independent of the choice
of contact form, and that the integrand involves only the second
and lower-order
derivatives (relative to a coframe field) while the lowest order
of local invariants is of order 4 as indicated by the Cartan tensor
$Q_{11}$.

Next, we will discuss the moduli space of $CR$-structures. Let
${\cal J}_{\xi}$ denote the space of all $CR$-structures compatible
with $\xi$. Let ${\cal C}_{\xi}$ denote the group of contact
diffeomorphisms with respect to $\xi$. Clearly, ${\cal C}_{\xi}$
acts on ${\cal J}_{\xi}$ by pulling back. The invariant $\mu_{\xi}$
is actually defined on the moduli space ${\cal J}_{\xi}/{\cal C}_{\xi}$.

Given a $CR$-structure $J$ in ${\cal J}_{\xi}$, we call a ``submanifold''
$S$ passing through $J$ a local slice if it is transverse to the orbit
of ${\cal C}_{\xi}$-action, so that any element in ${\cal J}_{\xi}$
near $J$ can be pulled back to an element of $S$ by means of a certain
contact diffeomorphism. In the early 1990's, Jack Lee and the author proved
the following:

\begin{indention}{2ex}
\noindent
\normalsize
\\
Local slices always exist for all cases (Cheng and Lee, 1995).  
\\
\end{indention}

As a corollary, the standard spherical $CR$-structure $[{\hat J}]$ in
${\cal J}_{\hat{\xi}}/{\cal C}_{\hat{\xi}}$ for $S^3$ is a strict
local minimum for $\mu_{\hat{\xi}}$ (Cheng and Lee, 1995).

Let $Q_{J}=2Re[i{Q_1}^{\bar 1}{\theta^1}{\otimes}Z_{\bar 1}]$. It is 
a straightforward computation to obtain the first variation formula:
${\delta}{\mu_{\xi}}(J)=-\frac{1}{8{\pi^2}}Q_{J}$. Consider the
downward gradient flow for ${\mu_{\xi}}$:

\begin{equation}
{\partial_t}J_{(t)}=Q_{J_{(t)}}.
\end{equation}

Since ${\delta}Q_J$ is subelliptic modulo the action of our symmetry
group ${\cal C}_{\xi}$, we can play a suitable ``De-Turck trick'' to
break the symmetry and imitate the usual $L^2$-theory for elliptic
operators to obtain the short time solution of Eq.(4)(Cheng and Lee,
1990). 
However, we can not prove the long term solution and convergence
even for $M=S^{3}$. This is related to the so-called Smale
conjecture as first pointed out by Eliashberg.

The Smale conjecture
asserts that the diffeomorphism group of $S^3$ is homotopy-equivalent
to the orthogonal group $O(4)$. Suppose we have the long term solution
and convergence of Eq.(4) for $M=S^{3}$. Then, any starting $J$ must
converge to ${\hat J}$, the unique spherical $CR$-structure on $S^3$
(up to symmetry). Therefore, the (certain marked) $CR$ moduli space
${\cal J}_{\hat{\xi}}^{'}/{\cal C}_{\hat{\xi}}^{'}$ is
contractible. But ${\cal J}_{\hat{\xi}}^{'}$ is contractible, too.
It follows that ${\cal C}_{\hat{\xi}}^{'}$ is contractible. Then, with the
aid of contact geometry, we can confirm the Smale conjecture.

To learn more analytic techniques which can be used to tackle Eq.(4),
we have been working
on some comparatively easier flows like the $CR$ Calabi flow and the $CR$
Yamabe flow. For the $CR$ Yamabe flow, S.-C. Chang and the author deformed
a contact form in the 
direction of the Tanaka-Webster curvature:

\begin{equation}
{\partial_t}{\theta_{(t)}}={\cal W}{\theta_{(t)}}.
\end{equation}

In their present work, Chang and Cheng obtain a
Harnack estimate and (possibly) the long term solution for Eq.(5).
\bigskip

\section{The Moduli Space of Spherical $CR$-Structures}

Let ${\cal S}_{\xi}$ denote the space of all spherical $CR$-structures
compatible with $\xi$. Since the linearization of the Cartan tensor
is subelliptic modulo the action of ${\cal C}_{\xi}$, the virtual dimension
of ${\cal S}_{\xi}/
{\cal C}_{\xi}$: the moduli space of spherical $CR$-structures is
finite. Let $M$ be a circle bundle over
a closed surface of genus $g>1$ with the Euler class $e(M)<0$. Let
$Pic(g,c_{1})$, the universal Picard variety, denote the space of
all pairs $(L,N)$ in which $L$ is
a holomorphic line bundle over a Riemann surface $N$ of genus $g>1$ with
$c_{1}(L)=e(M)$ modulo an equivalence relation defined by diffeomorphisms.
In 1996 and 1997, I-Hsun Tsai and the author studied the relation between
${\cal S}_{\xi}/{\cal C}_{\xi}$ and $Pic(g,c_{1})$. We found the following:

\begin{indention}{2ex}
\noindent
\normalsize
\\
For an above-mentioned circle
bundle $M$, there is a diffeomorphism between 
${\cal S}_{\xi}/^{'}{\cal C}_{\xi}$ and $Pic(g,c_{1})'$.
(The prime means a suitably modified version.)
Moreover, $Pic(g,c_{1})'$ is a complex manifold of dimension
$4g-3$ (Cheng and Tsai, 2000).
\\
\end{indention}

Our above result is similar to describing a Teichmuller
space by means of conformal classes. It is known in Teichmuller theory that
we can pick up a unique hyperbolic metric as a representative for 
each conformal class. A similar situation occurs for our spherical
$CR$ manifolds. In fact, our theory for the universal Picard variety
has counterparts in Teichmuller theory as shown in Table 1.

\bigskip
\noindent {\bf {Table 1. Comparison of two theories}}
\begin{center}
\begin{tabular}{l|l}
Teichmuller space & universal Picard variety \\ \hline
conformal classes & spherical $CR$ circle bundles \\ \hline
Riemannian hyperbolic metrics & pseudohermitian hyperbolic geometries \\
\end{tabular}
\end{center}

\bigskip
\noindent $\bf{Local\:Rigidity\:of\:Spherical\:CR-Structures}$

\noindent ${\bf (Discrete\:Moduli:} dim{\cal S}_{\xi}/{\cal C}_{\xi}=0)$

\bigskip
Let $Aut_{CR}(S^{3})$ denote the $CR$-automorphism group of
$(S^{3},{\hat{\xi}},{\hat{J}})$,
which is known to be isomorphic to $SU(2,1)/center$. Let $\Gamma$
denote a fixed point free finite subgroup of $Aut_{CR}(S^{3})$. Then,
${\Gamma}{\backslash}S^{3}$ inherits both contact and (spherical)
$CR$-structures from $(S^{3},{\hat{\xi}},{\hat{J}})$. This induced
spherical $CR$-structure on ${\Gamma}{\backslash}S^{3}$ is locally
rigid; i.e. it has no nontrivial deformation. (The algebraic reason
is that $H^{1}({\Gamma},su(2,1))=0$, in which the group cohomology
has coefficients in the holonomy representation: developing
map composed with the adjoint representation) (Burns and Shnider, 1976).
On the other
hand, note that ${\Gamma}{\backslash}S^{3}$ has positive constant
Tanaka-Webster curvature and zero torsion. Now, generalizing using
an analytical method, we obtain the following:

\begin{indention}{2ex}
\noindent
\normalsize
\\
Let $(M,J)$ be a closed spherical
$CR$ three-manifold. Suppose there is a contact form such that
the torsion $A_{11}=0$ and ${\cal W}>0$, $4{\cal W}(5{\cal W}^{2}+
3{\Delta}_{b}{\cal W})-3|{\nabla_b}{\cal W}|_{\theta}^{2}>0$. Then,
$J$ is locally rigid (Cheng, 1999).
\\
\end{indention}

Next, we want to compare two ${\Gamma}{\backslash}S^{3}$.
Suppose ${\Gamma_1}{\backslash}S^{3}$ and ${\Gamma_2}{\backslash}S^{3}$
are diffeomorphic. How can we distinguish one spherical $CR$-structure
from the other one? (They have the same $\mu_{\xi}$-value.) To deal with
this problem, we borrow ideas from quantum physics. If we view $\mu_{\xi}$
as a Lagrangian (action, more accurately) in $2+1$ dimensions, spherical
$CR$-structures are just classical fields. Therefore, ``quantum
fluctuations'' should give us refined invariants. In practice, we compute
the partition function heuristically:

\begin{eqnarray}
{\cal Z}_{k}&=&{\int_{{\cal J}_{\xi}/{\cal C}_{\xi}}}{\cal D}[J]
              e^{ik{\mu_{\xi}}([J])}\nonumber \\
 &=& k^{-\frac{dim}{2}}({\cal Z}_{sc}+O(k^{-1}))\:(k\:large),\nonumber
\end{eqnarray}

\noindent in which ${\cal Z}_{sc}$ is called the semi-classical
approximation. Note that only classical fields make contributions to
${\cal Z}_{sc}$. By imitating the finite dimensional case, we can compute
the modulus of ${\cal Z}_{sc}$ (Cheng, 1995):

\begin{eqnarray}
|{\cal Z}_{sc}| & = & lim_{k{\rightarrow}{\infty}}k^{\frac{dim}{2}}
                           |{\cal Z}_{k}|\nonumber \\
 & = & {\Sigma}_{J:spherical}\left| \frac{det{\Box_J}}{det'{\delta}Q_{J}}
                             \right| ^{\frac{1}{2}},\nonumber
\end{eqnarray}

\noindent in which ${\Box_J}$ is a fourth-order subelliptic self-adjoint
operator related to the ${\cal C}_{\xi}$-action, and ${\delta}Q_{J}$, the
second variation of ${\mu_{\xi}}$, is also a fourth-order subelliptic
self-adjoint operator modulo the ${\cal C}_{\xi}$-action. We can regularize
two determinants via zeta functions. ($det'$ means taking a regularized
determinant under a certain gauge-fixing condition.) 

\bigskip
\noindent $\bf{Conjecture}$: If $J$ is spherical,
$$ Tor(J){\stackrel{def}{=}}\left| \frac{det{\Box_J}}{det'{\delta}Q_{J}}
                             \right| ^{\frac{1}{2}}
$$
\noindent is independent of any choice of contact form, i.e., 
a $CR$ invariant. 

\bigskip

We expect to use $Tor(J)$ to distinguish among ``contact lens'' (or
``$CR$ lens'')
spaces $\{
{\Gamma}{\backslash}S^{3}\}$. Also, we note that $Tor(J)$ is a
contact-analogue of Ray-Singer's analytic torsion while no
contact-analogue is known for the Reidemeister torsion.
\bigskip

\section{Counting the Number of Complex Structures}

This is another ``quantum level'' problem in our ongoing project.
We will discuss the problem of counting the number of complex
structures on a closed (compact without boundary) four-manifold.
We hope to view this number as the partition function of a certain
3+1 quantum field theory (QFT in short).
  
Let us begin with a 0+1 theory, i.e., a particle moving in a closed
manifold $N$. The Hamiltonian of such a theory with supersymmetry
is the Laplace-Beltrami operator $\Delta$. All quantum ground states
or vacua are cohomology classes of $N$, represented by harmonic
forms (=zero eigenforms of $\Delta$). Now suppose $f$ is a Morse
function on $N$. Consider ${\Delta}_{tf}$ in which $d$ is replaced
by $e^{-tf}de^{tf}$. When $t{\rightarrow}{\infty}$, the harmonic forms
of ${\Delta}_{tf}$ are concentrated near the critical points of f.
These are the classical ground states (Witten, 1982).

The harmonic form corresponding to a critical point $P$ has a small
correction due to another critical point $Q$ via the trajectories of
${\nabla}f$ from $P$ to $Q$. This is quantum mechanical tunnelling,
which describes the probability of the transition $P{\rightarrow}Q$.
The boundary operator of Witten's chain complex (See Witten (1982) or 
Atiyah (1988) for
a clear explanation.) is interpreted in terms of such tunnelling.
(The homology of Witten's chain complex can be shown to identify
with the homology of $N$.) Witten's idea was later adopted by Floer (1989) 
and applied to the infinite-dimensional case of the manifold of connections.

\medskip

Next, we will give a brief introduction to the Donaldson-Floer theory. It
is a 3+1 QFT. A ``field'' when restricted to the three-space $M$
in this theory is a connection (or gauge field) of a certain, say,
$SU(2)$ bundle over $M$. The Morse function as mentioned above is
the Chern-Simons functional defined on the space of connections in 
this case. The critical points consist of flat connections which are
the classical ground states. Through consideration of the
associated Witten complex, we obtain the so-called Floer homology or
cohomology group $HF(M)$. This is the space of quantum ground states
or vacua for this theory. Now, suppose we decompose a closed 4-manifold
$X$ along $M$ (say, a homology 3-sphere) as shown in Fig.1.

\vspace{0.5cm}
\centerline{
\unitlength 1.00mm
\linethickness{0.2mm}
\begin{picture}(20,20)(30,137)
\put(15.33,140.01){\oval(12.67,13.33)[lt]}
\put(15.33,139.34){\oval(12.67,13.33)[lb]}
\put(54.67,140.34){\oval(12.67,13.33)[rt]}
\put(54.67,139.67){\oval(12.67,13.33)[rb]}
\bezier{160}(15.67,132.67)(34.00,139.01)(54.00,133.01)
\bezier{160}(15.38,146.71)(34.13,143.90)(54.61,146.85)
\bezier{48}(33.86,145.20)(29.79,139.40)(33.74,135.94)
\bezier{44}(33.86,145.20)(36.82,140.14)(33.49,135.94)
\put(15,148){\em $X^{+}$}
\put(32,147){\em $X$}
\put(56,148){\em $X^{-}$}
\put(31,132){\em $M$}
\put(3,125){\normalsize \bf Fig.1. Decomposing X along M}
\end{picture}
}
\vspace{1.5cm}
\bigskip
\noindent where $X=X^{+}{\cup}_{M}X^{-}$. Let ${\Sigma}^{+}$(${\Sigma}^
{-}$, respectively) denote the set of restrictions on $M$ of all instantons
on $X^{+}$($X^{-}$, respectively). Then,
${\Sigma}^{+}$, ${\Sigma}^{-}$ form cycles
in $HF(M)$. The intersection number represents the algebraic number of
instantons on $X$, (assuming it is finite) the Donaldson invariant,
denoted as $Z(X)$. We can write 

$$ Z(X)=<vac(X^{+})|{\;}vac(X^{-})>,$$

\noindent in which the vacuum $vac(X^{+})=[{\Sigma}^{+}]$ and the
vacuum $vac(X^{-})=[{\Sigma}^{-}]$ are both elements of $HF(M)$.
Also $<{}|{}>$ denotes the middle-dmension intersection number.
In Witten (1988), Witten presented a Lagrangian for this theory so that
$Z(X)$ identifies with its partition function.

\medskip
\indent
Now, we can describe our $3+1\;QFT$. We put an auxiliary contact 
structure $\xi$ on our closed oriented three-manifold $M$. A ``field'' is
a complex structure with the restriction on $M$ being a $CR$-structure
compatible with $\xi$. Our Morse function is the ${\mu}_{\xi}$ which
we introduce in $\S$3.
Spherical $CR$-structures which are critical points of ${\mu}_{\xi}$
are our classical ground states in this theory. 

Let ${\Sigma}^{+}$(${\Sigma}^{-}$, respectively) denote the set of all
$CR$-structures compatible with $\xi$ on $M$, which can be extended
to a complex structure on $X^{+}$$(X^{-}$, respectively). Now,
what is
the
associated ``Floer'' homology group $HF(M,{\xi})$, i.e., the space
of quantum vacua, for this theory?  Since the
Hessian ${\delta^2}{\mu}_{\xi}$ at a spherical $J$ is subelliptic
modulo ${\cal C}_{\xi}$, the dimension of its negative eigenspace is
finite. Therefore, the Morse index is well defined. (We do not need the
relative Morse index as in the case of the Donaldson-Floer theory.)
As usual, ${\Sigma}^{\pm}$ form cycles $[{\Sigma}^{\pm}]$ in $HF(M,{\xi})$
by pushing along the gradient flow of ${\mu}_{\xi}$ and seeing which
critical points they ``hang'' on (Atiyah, 1988). The vacuum
$vac(X^{+})$($vac(X^{-})$,
respectively) is defined as the homology class $[{\Sigma}^{+}]$
($[{\Sigma}^{-}]$, respectively) in $HF(M,{\xi})$. 
Moreover, we define the quantity $Z_{\xi}(X)$ as

\begin{eqnarray}
 Z_{\xi}(X)&{\stackrel{def}{=}}&<vac(X^{+})|{}vac(X^{-})>\nonumber\\
     &{\stackrel{def}{=}}&intersection{\;}number{\;}of{\;}
[{\Sigma}^{+}]{\;}and{\;}[{\Sigma}^{-}].\nonumber
\end{eqnarray}

\noindent The sum of $Z_{\xi}(X)$ over the isomorphism classes of tight
contact structures, denoted as $Z(X)$, can be interpreted as the (algebraic)
number of complex
structures on $X$. We propose the following ``physical'' problem:

\bigskip
\noindent $\bf{Problem\:1}$. Find a Lagrangian for the above theory
so that its partition function identifies with $Z(X)$.

\bigskip
There are topological obstructions for $M$ to admit 
spherical $CR$-structures (Goldman, 1983). For instance, the three-torus
$T^3$ does not admit any spherical $CR$-structure (compatible with
any given contact structure $\xi$). Therefore, $HF_{\star}(T^{3},{\xi})=0$
for any $\xi$, and we can propose the following problem for ``nonexistence'':

\bigskip
\noindent $\bf{Problem\:2}$. Suppose $X=X^{+}{\cup}_{T^3}X^{-}$. Find
conditions on $X$ and, perhaps, $X^{\pm}$ such that $Z(X)=0$.

\bigskip
We still need to investigate the relation between $Z(X)=0$ and the 
nonexistence of complex structures. Another situation occurs when $M$
is the standard contact 3-sphere $(S^{3},{\hat{\xi}})$. This admits
only one compatible spherical $CR$-structure, namely, the standard one
$\hat{J}$, which is a strict local minimum for $\mu_{\hat{\xi}}$
modulo symmetry as mentioned in Section III. It follows that
$HF_{0}(S^{3},{\hat{\xi}})=Z$ and $HF_{k}(S^{3},{\hat{\xi}})=0$ for
$k{\neq}0$. Therefore, we can propose the following problem concerning
``global rigidity'':

\bigskip
\noindent $\bf{Problem\:3}$. Suppose $X=X^{+}{\cup}_{S^3}X^{-}$. Find
conditions on $X$ and, perhaps, $X^{\pm}$ such that $Z(X)=1$.

\bigskip
Note that any tight contact structure on $S^3$ is
isotopy-equivalent to ${\hat{\xi}}$ according to Eliashberg (1992).
Therefore, $Z(X)$ in Problem 3 is just $Z_{\hat{\xi}}(X)$.


\bigskip
\section{Monopoles and Contact Structures}

Recently, Kronheimer and Mrowka (1997) studied contact structures
on 3-manifolds
via the 4-dimensional Seiberg-Witten monopole theory. Here, we will outline
another approach by Cheng and Chiu (1999).

Given a contact 3-manifold $(M,{\xi})$ and a background
pseudohermitian structure $(J,{\theta})$, we can discuss a canonical
$spin^c$-structure $c_{\xi}$ on ${\xi}^{\star}$.
With respect to $c_{\xi}$, we will consider the equations
for our ``monopole'' $\Phi$ coupled to the ``gauge field''
$A$. Here, $A$, the $spin^c$-connection, is required to be
compatible with the pseudohermitian connection on $M$. The
Dirac operator $D_{\xi}$ relative to $A$ is identified with
a certain boundary $\bar {\partial}$-operator 
${\sqrt 2}({\bar {\partial}}_b^{a}+
({\bar {\partial}}_b^{a})^{\star})$. In terms of the
components $({\alpha},{\beta})$ of $\Phi$, our equations read as

\begin{equation}
\left\{ \begin{array}{c}
       ({\bar {\partial}}_b^{a}+
({\bar {\partial}}_b^{a})^{\star})({\alpha}+{\beta})=0\\
(or\:{\alpha}_{,{\bar 1}}^{a}=0,\:{\beta}_{{\bar 1},1}^{a}=0)\\
da(e_{1},e_{2})-{\cal W}=|{\alpha}|^{2}-|{\beta}_{\bar 1}|^{2},
\end{array} \right.
\end{equation}

\noindent where $A=A_{can}+iaI$ and $\cal W$ denotes the Tanaka-Webster
curvature. Our first step in understanding Eq.(6) is as
follows:

\begin{indention}{2ex}
\noindent
\normalsize
\\
Suppose the torsion $A_{11}=0$.
Also, suppose $\xi$ is symplectically
semifillable, and that the Euler class $e({\xi})$ is not a torsion class. 
Then, Eq.(6) has nontrivial solutions (i.e., $\alpha$ and $\beta$
are not identically zero simultaneously)(Cheng and Chiu, 1999).
\\
\end{indention}

On the other hand, the Weitzenbock-type formula gives a nonexistence
result for ${\cal W}>0$. Together with the above existence result,
we can conclude the following:

\begin{indention}{2ex}
\noindent
\normalsize
\\
Suppose the torsion $A_{11}=0$ and
the Tanaka-Webster curvature ${\cal W}>0$. Then, either $\xi$ is not
symplectically semifillable, or $e({\xi})$ is a torsion class (Cheng
and Chiu, 1999).
\\
\end{indention}

We note that Rumin (1994) proved that $M$
must be a rational homology sphere under the conditions given
above using a different method. Also, we do not know how
to deal with the solution space of Eq.(6) in general although
we hope that further study of Eq.(6) will produce invariants of contact
structures.

\bigskip

\noindent $\bf{Acknowledgment}$

\bigskip

In preparing this article the author benefited from a number of
conversations with I-Hsun Tsai and Chin-Lung Wang. Also, the author 
would like to thank Professors Chuu-Lian Terng and I-Hsun Tsai for
inviting him to lecture at the NCTS conference. The main part of this 
article is based on the author's notes for that talk. 
The research was partially supported by the National Science Council under
grant NSC 88-2115-M-001-015 (R.O.C.).

\bigskip

\noindent {\bf References}

\bigskip

Altschuler, S. J. (1995)
A geometric heat flow for one-forms on three-dimensional
manifolds.
{\em Illinois J. of Math.}, {\bf 39}, 98-118.

Atiyah, M. F. (1988)
New invariants of three and four dimensional manifolds.
{\em The Mathematical Heritage of Herman Weyl, Proc.
Symp. Pure Math.}, {\bf 48}, 285-299.

Burns, D. and C. Epstein (1988)
A global invariant for three dimensional CR-manifolds.
{\em Invent. Math.}, {\bf 92}, 333-348.

Burns, D. and C. Epstein (1990)
Characteristic numbers of bounded domains.
{\em Acta Math.}, {\bf 164}, 29-71.

Burns, D. and S. Shnider (1976)
Spherical hypersurfaces in complex manifolds.
{\em Invent. Math.}, {\bf 33}, 223-246.

Cartan, E. (1932a)
Sur la g$\acute{e}$ometrie pseudo-conforme des hypersurfaces de
l'espace de deux variables complexe I.
{\em Ann. Mat.}, {\bf 11}, 17-90.

Cartan, E. (1932b)
Sur la g$\acute{e}$ometrie pseudo-conforme des hypersurfaces de
l'espace de deux variables complexe II. 
{\em Ann. Sc. Norm. Sup. Pisa}, {\bf 1}, 333-354.

Cheng, J. H. and H. L. Chiu (1999)
{\em Monopoles and contact 3-manifolds}.
Part of NSC Report for Project NSC 88-2115-M-001-015,
National Science Council, R.O.C., Taipei, Taiwan, R.O.C..

Cheng, J. H. (1987)
On the curvature of CR-structures and its covariant derivatives.
{\em Math. Z.}, {\bf 196}, 203-209.

Cheng, J. H. (1999)
Rigidity of automorphisms and spherical CR structures.
{\em Proc. A.M.S.}, (to appear).

Cheng, J. H. (1995)
The geometry of Cauchy-Riemann manifolds in dimension three.
{\em Proceedings of the Second Asian Mathematics Conference},
Nakhon Ratchasima, Thailand.

Cheng, J. H. and J. M. Lee (1990)
The Burns-Epstein invariant and deformation of CR structures.
{\em Duke Math. J.}, {\bf 60}, 221-254.

Cheng, J. H. and J. M. Lee (1995)
A local slice theorem for 3-dimensional CR structures.
{\em Amer. J. Math.}, {\bf 117}, 1249-1298.

Chen, B. Y. and H. S. Lue (1981)
Chern classes and Bochner-Kaehler metrics.
{\em Coll. Math.}, {\bf 44}, 203-208.

Chern, S. S. and J. K. Moser (1974)
Real hypersurfaces in complex manifolds.
{\em Acta Math.}, {\bf 133}, 219-271.

Chern, S. S. and J. K. Moser (1983)
Erratum: Real hypersurfaces in complex manifolds.
{\em Acta Math.}, {\bf 150}, 297.

Cheng, J. H. and I-H. Tsai (2000)
Deformation of spherical CR structures and the universal
Picard variety.
{\em Commun. Anal. and Geom.}, (to appear)

Eliashberg, Y. (1992)
Contact 3-manifolds twenty years since J. Martinet's work.
{\em Ann. Inst. Fourier}, {\bf 42}, 165-192.

Fefferman, C. (1974)
The Bergman kernel and biholomorphic mappings of pseudoconvex
domains.
{\em Invent. Math.}, {\bf 26}, 1-65.

Floer, A. (1989)
An instanton invariant for 3-manifolds.
{\em Commun. in Math. Phys.}, {\bf 118}, 215-240.

Goldman, W. (1983)
Conformally flat manifolds with nilpotent holonomy and 
the uniformization problem for 3-manifolds.
{\em Trans. Amer. Math. Soc.}, {\bf 278}, 573-583.

Gray, J. W. (1959)
Some global properties of contact structures.
{\em Ann. Math.}, {\bf 69}, 421-450.

Hamilton, R. S. (1982)
The inverse function theorem of Nash and Moser.
{\em Bull. Amer. Math. Soc.}, {\bf 7}, 65-222.

Kronheimer, P. B. and T. S. Mrowka (1997)
Monopoles and contact structures.
{\em Invent. Math.}, {\bf 130}, 209-255.


Lee, J. M. (1986)
The Fefferman metric and pseudohermitian invariants.
{\em Trans. Amer. Math. Soc.}, {\bf 296}, 411-429.

Lutz, R. (1977)
Structures de contact sur les fibre's principaux en cercles de
dimension 3.
{\em Ann. Inst. Fourier}, {\bf 3}, 1-15.

Martinet, J. (1971)
Formes de contact sur les vari$\acute{e}$ti$\acute{e}$s de 
dimension 3.
{\em Lect. Notes in Math.}, {\bf 209}, 142-163.

Rumin, M. (1994) 
Formes diff{\'e}rentielles sur les vari{\'e}t{\'e}s
de contact.
{\em J. Diff. Geom.}, {\bf 39}, 281-330.

Tanaka, N. (1975)
{\em A Differential Geometric Study on Strongly
Pseudo-Convex Manifolds},
Kinokuniya Co. Ltd., Tokyo

Webster, S. M. (1977)
On the pseudo-conformal geometry of a K$\ddot{a}$hler manifold.
{\em Math. Z.}, {\bf 157}, 265-270.

Webster, S. M. (1978)
Pseudohermitian structures on a real hypersurface.
{\em J. Diff. Geom.}, {\bf 13}, 25-41.

Witten, E. (1982)
Supersymmetry and Morse theory.
{\em J. Diff. Geom.}, {\bf 17}, 661-692.

Witten, E. (1988)
Topological quantum field theory.
{\em Commun. in Math. Phys.}, {\bf 117}, 353-386.
}
\end{document}